\documentclass[graybox]{svmult}

\usepackage{mathptmx}       
\usepackage{helvet}         
\usepackage{courier}        
\usepackage{makeidx}         
\usepackage{graphicx}        
\usepackage{multicol}        

\makeindex             
                       
\usepackage[active]{srcltx}
\usepackage{amsmath,amssymb}
\usepackage{graphicx}
\usepackage{float}
\newcommand{\bb}[1]{{\mathbb #1}}
\newcommand{\mc}[1]{{\mathcal #1}}
\def\T{\bb T}

\def\TmN{T^N_{\rm max}}
\newcommand{\p}{\partial}

\begin{document}

\title*{Interacting particle systems: hydrodynamic limit  versus high density limit}
\author{Tertuliano Franco}
\institute{Tertuliano Franco \at Universidade Federal da Bahia, Salvador, Brazil, \email{tertu@impa.br}}

%
%
\maketitle

\abstract*{In this short survey we compare aspects of two different approaches for scaling limits of interacting particle systems,  the \emph{hydrodynamic limit} and the \emph{high density limit}.  We present some examples, comments and open problems on each approach for different scenarios: the law of large numbers,  the central limit theorem, and  the large deviations principle. It is given a special attention to a recent result \cite{FrancoGroisman} about the high density limit of a system exhibiting explosion of the total numbers of particles in finite time. }

\abstract{In this short survey we compare aspects of two different approaches for scaling limits of interacting particle systems,  the \emph{hydrodynamic limit} and the \emph{high density limit}.  We present some examples, comments and open problems on each approach for different scenarios: the law of large numbers,  the central limit theorem, and  the large deviations principle. It is given a special attention to a recent result \cite{FrancoGroisman} about the high density limit of a system exhibiting explosion of the total numbers of particles in finite time. }

\section{Introduction}
\label{sec:1}

A central question in Statistical Mechanics is about the passage from discrete systems to the \emph{continuum}. And, consequently,  how the intrinsic properties of a discrete systems are inherited by the \emph{continuum}. Aiming for rigorous results in this scope, fruitful mathematical theories have been developed since the last century. 

An important class of discrete systems are the so-called \emph{interacting particle systems}. Roughly speaking, an interacting particle system is a discrete system that evolves in time according to random clocks under some interaction among particles. To clarify ideas, two  examples of interacting particle systems are presented in the Section \ref{sec:2}.

A quantity of interest associated to a particle system is the spatial density of particles. Since the system evolves in time, its spatial density of particles evolves as well. Is therefore natural to investigate the possible limits for the \emph{time trajectory} of the spatial density of particles.

The limiting object for the time trajectory of the spatial density of particles is usually described as the solution of some partial differential equation. A standard hypothesis is to suppose that, at initial time, the spatial density of particles converges to a profile $\varphi$ as the mesh of the lattice goes to zero. This profile $\varphi$ will be, as reasonable, the initial condition of the respective partial differential equation.

The nature of the convergence (topology, parameters to be rescaled) is the subject of this paper. For sake of clarity, we take the liberty to divide the main types of convergence in two classes: the hydrodynamic limit and the high density limit. The expression \emph{hydrodynamic limit} is widely used in the literature. On the hand, the \emph{high density limit} nomenclature is less known, being employed in the paper \cite{Kotelenez88}. 

The hydrodynamic limit  consists of the limit for the time trajectory of the 
spatial density of particles where the parameters to be rescaled are time and space. The space among particles is lead to zero, while the time the system has to evolve is lead to infinity. If the time is taken as the inverse of the space among particles, this situation is called the \emph{ballistic scaling}. If the time is as taken as square of the inverse of the space among particles, it is called the \emph{diffusive scaling}. The initial configuration of particles is randomly chosen according to a distribution related to the fixed profile $\varphi$. 

The high density limit consists of the limit for the time trajectory of the 
spatial density of particles where the parameters to be rescaled are: time, space \textit{and initial quantity of particles per site}. While space among particles is lead to zero,  time and initial quantity of particles per site are lead to infinity. As suggested, the fact the initial quantity of particles increases in a meaningful way originated the nomenclature.

A vast literature has been produced about the hydrodynamic limit, which is nowadays an exciting and active research area. For a reference in the subject, we refer the reader to the classical book \cite{kl}. For very important techniques in  the area, we cite the \emph{Entropy Method}, the \emph{Relative Entropy Method}, non-gradient techniques, attractiveness techniques, among many others. In its turn, the high density limit approach had important papers about as \cite{at, Blount91,Blount92,Blount96,Kotelenez86,Kotelenez88}. As a more recent paper on the subject, we cite \cite{FrancoGroisman}.  

The hydrodynamic limit is far more studied and known. But the high density limit has interesting characteristics and a plenty of open problems, some of them discussed here. Our goal in this short survey is to compare main aspects of each approach, exemplify them, and debate for whose models each one is more suitable to (in some sense).

The outline of this paper is: in Section \ref{sec:2}, two interacting particle systems are presented. In Section \ref{sec:3},  it is made a  comparison of results in each approach for the law of large numbers scenario. In Section \ref{sec:4}, the same for the central limit theorem scenario, and in Section \ref{sec:5}, the same for the large deviations principle scenario.

\section{Two interacting particle systems}
\label{sec:2}

In this section we present two dynamics of interacting particle systems. For a classical reference on particle systems, we cite the book \cite{L.}. The first example is the symmetric simple exclusion process, the second one is a system of independent random walks with birth-and-death dynamics. Denote by 
$$T_N\;=\;\bb Z/ (N\bb Z)\;=\;\{0,1,2,\ldots, N-1\}$$
the discrete one-dimensional torus with $N$ sites.

\runinhead{The symmetric simple exclusion process}
The symmetric simple exclusion process, abbreviated by SSEP, is a quite standard, widely studied model in Probability and Statistical Mechanics. In words, in the SSEP each particle performs an independent  continuous-time symmetric random walk except when some particle tries to jump to an already occupied site. When this happens, this jumps is forbidden, and nothing happens. This  \emph{exclusion rule}  originates the name exclusion process. In several models of physical phenomena, fermions dynamics are constrained by an exclusion rule. 

Of course, since a particle can not jump to an already occupied site, the state space in this case is 
$\{0,1\}^{T_N}$. For $x\in T_N$, we will write down $\eta(x)$ for the number of particles at the site $x$ in the configuration $\eta$ of particles.  See Figure \ref{fig:3} below, where black balls represent particles.
\begin{figure}[H]
\sidecaption
\includegraphics[scale=0.65]{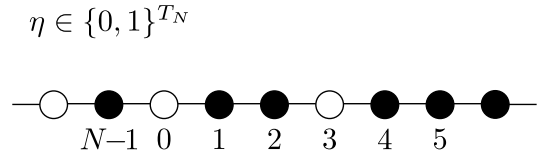}
\caption{A configuration  of particles $\eta\in\{0,1\}^{T_N}$. Observe that $\eta(N-1)=1$, $\eta(0)=0$, 
$\eta(1)=1$, $\eta(2)=1$, \emph{etcetera}. Notice that, since $T_N$ is the discrete torus,  $x=0$ and $x=N$ represent the same site.}
\label{fig:3}      
\end{figure}

\begin{figure}[H]
\sidecaption
\includegraphics[scale=0.65]{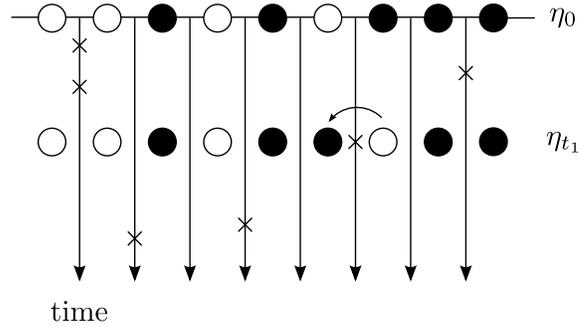}
\caption{At right, an evolution of the initial configuration $\eta_0$ according to the Poisson processes (the marks represent the time arrivals). At time $t_1$, a particle jumps to a neighbor site. Notice that at the three marks in times previous to $t_1$ nothing happened because both sites related to the mark were empty or occupied.}
\label{fig:2}      
\end{figure}

The dynamics is the following: to each \emph{edge} $(x,x+1)$ of the discrete torus $T_N$, it is associated a Poisson point process\footnote{A Poisson process can be described as marks in time, being the time between marks i.i.d of exponential distribution. The parameter $N^2$ has to do with the scaling we are going to perform later.} of parameter $N^2$, all of them independent. At a time arrival of the Poisson process corresponding to the edge $(x,x+1)$, the occupations at $\eta(x)$ and $\eta(x+1)$ are interchanged. In case that both sites $x$ and $x+1$ are occupied, of course, nothing happens\footnote{Corroborating the exclusion rule.}. Since the Poisson processes are independent, the probability to observe two marks at same time is zero. Hence, there is no chance to a particle be in ``doubt'' whether to jump, and the construction is well defined. Figure \ref{fig:2} illustrates ideas. 
Given an initial configuration of particles $\eta\in\{0,1\}^{T_N}$, this construction yields the continuous time Markov chain $\{\eta_t\;;\; t\geq 0\}$, which is the so-called SSEP. 

\runinhead{Independent random walks with birth-and-death dynamics} In this particle system we have a superposition of two standard dynamics. One is given by independent random walks, where each particle has its own random clock (a Poisson process associated to him). When the clock rings, the particle chooses with equal probability one of the neighbor sites and jumps to there. For independent random walks, there is no interaction among particles. The other part of the dynamics is given by birth and death of particles at each site. The  birth and death rates at a site $x$ are given by functions $b$ and $d$, respectively, of the number of particles at that site $x$.

In each site of $T_N$ we allow a nonnegative integer quantity  of indistinguishable particles. A configuration of particles will be denote by $\eta$, which is an element of $\{\mathbb{N}\cup \{0\}\}^{T_N}$. As before, we will write down $\eta(x)$ for the number of particles at the site $x$ in the configuration $\eta$ of particles, see Figure \ref{fig:1} below.
\begin{figure}[H]
\sidecaption
\includegraphics[scale=0.65]{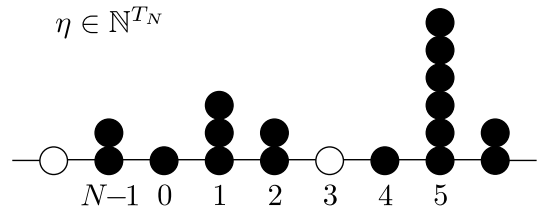}
\caption{A configuration of particles $\eta\in\mathbb{N}^{T_N}$. Observe that $\eta(N-1)=2$, $\eta(0)=1$, 
$\eta(1)=3$, $\eta(2)=2$ and $\eta(3)=0$. Again, since $T_N$ is the discrete torus,  $x=0$ and $x=N$ represent the same site.}
\label{fig:1}      
\end{figure}

Next, we construct the system of independent random walks with birth-and-death dynamics. 
Fix two nonnegative smooth\footnote{The smoothness is not necessary here. It will be required only in the later scaling.} functions $b,d:\bb R_+\to\bb R_+$ such that $d(0)=0$ and fix $\ell=\ell_N$ a positive parameter. In the next section we will see that this parameter $\ell$ represents the number of particles per site at the initial time.

Consider the following transition rates:
\begin{enumerate}
 \item[\textbullet] at rate $N^2\eta(x)$, a particle jumps from $x$ to $x+1$; 
 \item[\textbullet] at rate $N^2\eta(x)$, a particle jumps from $x$ to $x-1$; 
\item[\textbullet] at rate $\ell b(\ell^{-1}\eta(x))$, a new particle is created at $x$; 
\item[\textbullet] at rate $\ell d(\ell^{-1}\eta(x))$, a particle is destroyed at $x$. 
\end{enumerate}
The transitions above are assumed for all $x\in \T_N$. Each transition  corresponds to an arrival of an independent Poisson process of  parameter given by the respective rate. That is, an analogous graphical construction (as the aforementioned for the SSEP) can be made in this case, see \cite{FrancoGroisman}. 

 Since there are no assumptions about the growth of $b$, the waiting times of this Markov chain can be summable\footnote{Meaning that the total quantity of particles has exploded. For more on explosions of Markov chains see \cite{Norris}.}. In this case we say that the process \emph{explodes} or \emph{blows up}, and we define the state of the process as $\infty$ for times greater or equal than the sum of all the waiting times, that we call $\TmN$. More precisely, define
\[
 \tau^N_y:=\inf\{t\ge 0\colon \|\eta(t)\|_\infty\ge y \}\quad  \mbox{ and }\quad 
\tau^N_{\textrm{blow-up}}:=\lim_{y\to \infty}  \tau^N_y\,.
 \]
Then, for $t<\tau^N_{\textrm{blow-up}}$, we define $\eta(t)$ by means of the rates  stated before, and
for $t\ge \tau^N_{\textrm{blow-up}}$, we define $\eta(t)=\infty$. This characterizes a continuous time Markov chain 
$$\{\eta_t\;;\;t\geq 0\}$$
with state space $\bb N^{T_N}\cup \{\infty\}$.

\section{Law of large numbers scenario for each setting}\label{sec:3}

In order to state the limit for the time trajectory of the spatial density of particles, we need to define first what we mean by a spatial density of particles. We point out that the definition of the spatial density of particles is different for each setting, the hydrodynamic limit or the high density limit. In the first one, given the SSEP $\{\eta_t\;;\; t\geq 0\}$ described in the Section \ref{sec:2}, the spatial density of particles,  usually called the \emph{empirical measure}, is 
defined as 
\begin{equation}\label{emp_measure}
 \pi^N_t\;:=\;\frac{1}{N}\sum_{x\in T_N}\eta_{t}(x)\,\delta_{\frac{x}{N}}\,.
\end{equation}
As can be easily seen, the empirical measure is:
\begin{itemize}
 \item a positive measure (since it is a sum of deltas of Dirac);
  \item a random measure (since  it is a function of $\eta_t$, which is random);
 \item   a measure with total mass bounded by one (by the normalization constant $N$);
 \item a measure that gives mass $1/N$ at the point $x/N$ belonging at the continuous one-dimensional torus $$\bb T=\bb R/\bb Z=[0,1)$$ if there is a particle
 at the site $x\in T_N$ (at that time $t$), and gives measure $0$ otherwise.
 \end{itemize}
 Denote by $\mc M_+$ the space of positive measures on $\bb T$ with mass bounded by one and by $\mc D([0,T], \mc M_+)$ the set of c\`adl\`ag\footnote{From the French, right continuous and  with side left limits.} time trajectories taking values on $\mc M_+$. We notice that the time trajectory of the empirical measure \eqref{emp_measure} is a random element taking values in  $\mc D([0,T], \mc M_+)$.

The theorem stated next is what we call the hydrodynamic limit (for the SSEP). For a proof, see \cite[Chapter 4]{kl}. We notice that the topology assumed in the convergence in  distribution ahead is the Skorohod topology on $\mc D([0,T], \mc M_+)$. For an exposition on the Skohorod topology, see  \cite{Bill} or \cite{kl}.

\begin{svgraybox}
\begin{theorem} \label{th:hlrm}
Fix $\varphi: \mathbb{T} \to [0,1]$ a smooth function and $T>0$. Suppose that the initial distribution of particles for the SSEP are chosen in such a way, as $N\to\infty$,
\begin{equation*}
 \pi^N_0\longrightarrow \varphi(u)\,du\quad \textrm{in probability}. 
\end{equation*}
 Then,   as $N\to\infty$,
 \begin{equation*}
  \{\pi^N_{t}\;;\;0\leq t\leq T\}\longrightarrow \{\rho(t,u)\,du\;;\;0\leq t\leq T\}\quad \textrm{in distribution,}
 \end{equation*}
  where $\rho(t,u)$ is the  unique solution of the periodic heat equation with initial condition $\varphi$, or else,
 \begin{equation}\label{heat}
\left\{
\begin{array}{ll}
 \partial_t \rho(t,u) \; =\; \p_{uu} \rho(t,u)\,, \qquad& t \geq 0,\, u\in \bb T\,,\\
  \rho(0,u) \;=\; \varphi(u), &u \in \bb T\,.
\end{array}
\right.
\end{equation}

\end{theorem}
\end{svgraybox}
It is not our intention to resume a huge research area in few words. Nevertheless, let us make some remarks.  The hydrodynamic limit has been successfully applied in interacting particles systems as: the zero range model; the symmetric simple exclusion process, the asymmetric simple exclusion process, the Ginzburg-Landau model, the generalized exclusion process, among others. 

We point out here the power of the existing methods for proving the hydrodynamic limit of models whose microscopic interactions often lead to non-linear partial differential equations. In the above case, it is obtained a linear heat equation. For the zero range process, the partial differential equation would the a non-linear one, with $\p_{uu}\Phi(\rho)$ replacing $\p_{uu}\rho$ in \eqref{heat}, where the function $\Phi$ is defined via the microscopic interaction.

Despite the wide applicability, the available techniques for hydrodynamic limit are not suitable for systems with huge birth's rate of new particles. There are some papers on the subject as \cite{Mourragui96}, but in a scheme where the  birth's rate of new particles is small in some sense. The reason is simple, 
in the situation where the particle system explodes in finite
time, the expectation of the number of particles is infinity at any positive
time. Hence, any method based on expectations is doomed to fail. And most of hydrodynamic techniques are based on expectation techniques.

For $x\in T_N$, let $u_x=x/N\in \bb T$. We define now $X^N:\bb R_+\times \bb T\to\bb R_+$, the spatial density of particles   for the high density limit scenario, by
\begin{equation*}
 X^N(t,u_x) \;=\; \ell^{-1} \eta_t(x)\,.
\end{equation*}
For $u_x<u<u_{x+1}$, we define the spatial density via a linear interpolation, i.e.,
\begin{equation*}
 X^N(t,u) \;=\; (Nu-x)\; X^N(t,u_{x+1}) + (x+1-Nx)\; X^N(t,u_k)\,.
\end{equation*}
If $\eta(t)=\infty$, we say that $ X^N(t,\cdot )=\infty$ as well. 

Before stating a high density limit theorem, let us say some words about the partial differential equation
\begin{equation}
\label{PDE}
\left\{\begin{array}{ll}
\p_t\rho(t,u) = \p_{uu} \rho(t,u) + f(\rho(t,u)) \,, \qquad & t\in [0,T),\, u\in\bb T\,,\\
 \rho(0,u) = \varphi (u) \ge 0\,, \qquad & u \in \bb T\,. \\
\end{array}\right.
\end{equation}
Above, $f=b-d$, where $d$ and $b$ are the aforementioned smooth functions that drive the birth and death of particles, and $\varphi$ is smooth and nonnegative. Since there is not restriction about the growth of $f$, the solution of \eqref{PDE} can exhibit a phenomena called \emph{blow-up} or \emph{explosion}. In this case, there is a finite time $T_{\textrm{blow-up}}$ for which
\begin{equation}\label{eq03}
\lim_{t\nearrow T_{\textrm{blow-up}}} \|\rho(t,\cdot)\|_{\infty} = \infty\,, 
\end{equation}
and such that $\|\rho(t,\cdot)\|_{\infty}$ is finite for times smaller than  $T_{\textrm{blow-up}}$. In this case, the solution $\rho$ of \eqref{PDE} is defined only in the time interval $[0,T_{\textrm{blow-up}})$. If there is no explosion, we would say that $T_{\textrm{blow-up}}=\infty$. 

Remark: a well known condition on the nonlinear term $f$ that assures the existence of
solutions with blow-up is being convex, strictly positive in some interval $[a,
+\infty)$ and $\int^\infty_a \frac{ds}{f(s)} < \infty$.
\begin{svgraybox}
\begin{theorem}[F. , Groisman '12] Assume that
\begin{enumerate}
 \item[(A1)] $\lim_{N\to\infty}\Vert X^N(0,\cdot)-\varphi(\cdot)\Vert_\infty = 0$ almost
surely;
 \item[(A2)] for any  $c>0$,  $\sum_{N\geq 0} N^3
e^{-c\, \ell(N)}<\infty\,$.
\end{enumerate}
Then, for any $T\in [0,T_{\textrm{blow-up}})$,
\begin{equation*}\label{eq2}
\lim_{N\to\infty} \;\sup_{t\in[0,T]} \Vert  X^N(t,\cdot)-\rho(t,\cdot)\Vert_\infty\;=\; 0\quad \textrm{ almost surely,}
\end{equation*}
where $\rho$ is the solution of \eqref{PDE} and $T_{\textrm{blow-up}}$ is given in \eqref{eq03}.
\end{theorem}
\end{svgraybox}
The assumption (A1) allows to interpret  $\ell$ as the quantity of particles per site at initial time. Roughly speaking, since $X^N_0$ is $\ell^{-1} \eta_0$,  and $X^N_0$ converges in the supremum norm to the function $\varphi$, the initial quantity of particles at a point $u\in \bb T$ is 
 of order $ \ell \,\varphi(u)$. 
 
 We remark that $\pi^N_t$ and $X^N(t,u)$ are equivalent in some sense. With due care, the hydrodynamic limit could be stated in terms of piecewise affine functions and, \emph{vice versa}, the high density limit could be stated in terms of a suitable empirical measure. 
 
 The result above is proved in \cite{FrancoGroisman} making use of couplings. It is a challenging problem to obtain some analogous results in the hydrodynamic limit setting. For instance, the zero range process has no limit of particles per site. Superposing this dynamics with birth of particles, explosions can occur under suitable choice of rates. The hydrodynamic limit of this model could be studied.

 The known techniques for the high density limit are strongly support on three pillars:  martingales, Duhamel's Principle and smoothing properties of the heat equation semi-group. Duhamel's Principle is a general idea widely applied in ordinary different equations, partial differential equations, numerical schemes, \emph{etcetera}. For a system whose dynamics has two parts, being one of them linear, the solution (in time) can be expressed as the initial condition evolved by the linear part plus an convolution of the nonlinear part with the semi-group of the linear part.
 
 Keeping this in mind, we can realize why the literature on high density limit is concentrated in dynamics involving independent random walks: in order to apply Duhamel's Principle, it is necessary to have a linear part in the dynamics, and the independent random walks plays this role. It is an open problem to extend the high density limit for others dynamics as the zero range process, for example.
 
In the paper \cite{Mourragui96}, it was considered the exclusion process superposed with birth dynamics, but with no explosions.
 It is  a challenging problem to prove the hydrodynamic limit for the zero range process with birth of new particles and explosion in finite time.
\section{Central limit theorem  scenario for each setting}
\label{sec:4}

There are a lot of results on the central limit theorem for both settings, the hydrodynamic limit and the high density limit. In the high density limit setting, we cite \cite{Blount91,Blount96,Kotelenez86,Kotelenez88}. In the hydrodynamic limit setting, we cite \cite{g,jl}.

In both settings, the limit is usually described through generalized Ornstein-Uhlenbeck processes, see \cite{HS} or \cite{kl} on this type of stochastic process. 

We cite as an interesting open problem to prove the central limit theorem for the spatial density of particles near the explosion for the model considered in \cite{FrancoGroisman}.

\section{Large Deviation Principle  scenario for each setting}
\label{sec:5}
Recall that $\mc M_+$ denote the space of positive measures on $\bb T$ with mass bounded by one and by $\mc D([0,T], \mc M_+)$ the set of c\`adl\`ag time trajectories with values on $\mc M_+$. 

By a large deviations principle in the hydrodynamic setting, we mean the existence of a lower semicontinuous rate function 
$$I:\mc D([0,T], \mc M_+)\to \bb R_+\cup\{\infty\}$$
such that:
\begin{svgraybox}
 \it For each closed set $\mc C$, and each open set $\mc O$ of $\mc D ([0,T],\mc M_+)$, 
\begin{eqnarray*}
 \limsup_{N\to\infty} \frac{1}{N}\log Q^N[\mc C] \;\leq\; -\inf_{\pi\in \mc C} I(\pi)\,,\\
 \liminf_{N\to\infty} \frac{1}{N}\log Q^N[\mc O] \;\geq\; -\inf_{\pi\in \mc O} I(\pi)\,,
 \end{eqnarray*}

\end{svgraybox}
\noindent
where $Q^N$ is the probability induced the in the space $\mc D([0,T], \mc M_+)$ by the empirical measure. A proof of the large deviation principle for the SSEP can be found in \cite[Chapter 10]{kl}. 

For the high density limit, the statement of a LDP is analogous, \emph{mutatis mutandis} with respect to the topology.
However, there is no LDP available yet for the high density limit scenario and we list some difficulties for proving it. 

An important step in order to obtain a LDP for a model is to obtain a law of large numbers for a class of perturbed process. The Radon-Nykodim between the original process and the perturbed one will further give the cost to observe the limit given by the perturbed one from the point-of-view of the original process. In the proof of a LDP, it is made a careful analysis and precise optimization over the perturbations.

For the SSEP, the perturbations are given by weakly asymmetric exclusion process, where the asymmetry is driven by a smooth function $H$, being the limit for the hydrodynamic limit given by a solution of a heat equation with a non-linear Burgers term added, see \cite[page 273]{kl}.

As commented in the Section \ref{sec:3}, the known techniques for the high density limit are strongly based on martingales, Duhamel's Principle and smoothing properties of the semi-group corresponding to the linear part of the dynamics. In order to obtain a LDP for the high density limit (let us say, in the case of independent random walks with birth-and-death dynamics) it would be necessary to prove the high density limit for some non-linear situation not attained yet in the literature. 

Other obstacle is the presence of two superposed dynamics. To observe a given profile that differs from the expected limit, it would be possible to consider two different perturbations at same time, one about the diffusion and another about  the birth-and-death of particles. Performing variational analysis over two competing different perturbations is a complicated situation in LDP. 

In resume, LDP is a challenging open problem in the high density limit scenario.

\begin{acknowledgement}
T.F. was supported through a grant ``BOLSISTA DA CAPES - Bras\'ilia/Brasil" provided by CAPES (Brazil) and through a grant ``PRODOC-UFBA (2013-2014)''. The author also would like to thank the anonymous referees for valuable comments and improvements.
\end{acknowledgement}

\bibliographystyle{amsplain}

\end{document}